# Coherent Measures of Discrepancy, Uncertainty and Dependence, with Applications to Bayesian Predictive Experimental Design

A. P. Dawid
Department of Statistical Science
University College London

August 3, 1998

**SUMMARY**

We show how, associated with any decision problem, we may derive related functions measuring uncertainty, discrepancy and dependence of distributions. Such "coherent" functions have special properties, which we characterise, and each functions essentially determines the others. The theory is applied to the Bayesian formulation of the problem of choosing an experiment in order to make a subsequent prediction. It is shown that coherent choice criteria may be based on any of the coherent functions, with related functions yielding identical solutions.

# 1 Introduction

Let $Z$ be an uncertain quantity with values in $\mathcal{Z}$. We may regard $Z$ as an observable outcome of an experiment, or as the parameter of a statistical model. Using Bayesian decision theory, if $Z$ is assigned a distribution $P$, any decision problem, with loss function $L$ depending on the outcome of $Z$, may be solved.

In this paper we show how such a "terminal decision problem" also determines functions which may be used to quantify a variety of concepts relating to a space $\mathcal{P}$ of distributions for $Z$: a *scoring rule* for assessing the empirical adequacy of $P \in \mathcal{P}$; the *uncertainty* in a distribution $P \in \mathcal{P}$; the *discrepancy* between two distributions $P$ and $Q \in \mathcal{P}$; and a measure of the *dependence* of $Z$ on another variable. These functions inherit various special properties, which we here characterise, from their decision-theoretic origins.

The Bayesian approach to decision theory can also be applied to solve the problem of choosing an experiment to perform, before proceeding to the terminal decision problem. It leads to a criterion for such choice based on the minimisation of overall expected loss. We show how this procedure can be re-expressed in terms of any of the derived functions, leading in each case to a naturally appealling criterion. The basic structure of the design problem considered here is laid out in §2. It is formulated very generally, in terms of a terminal decision problem involving prediction of a future observable $Z$; the more common problem of estimation can be tackled as a special case of this.

In §3, we introduce our general formulation of a terminal decision problem, governed by the outcome of $Z$, and in particular show how it gives rise to a measure of uncertainty in a distribution. Section 4 then applies this to the problem of experimental design.

In §5 we describe how a *proper scoring rule* for assessing a probability distribution for $Z$ may be associated with our terminal decision problem, and how this may again be used as the basis of the design problem. Section 6 derives from this a measure of discrepancy between two distributions for $Z$, demonstrating how the design criterion may again be re-expressed in terms of this discrepancy function. In Section 7 we introduce a measure of the dependence of $Z$ on another quantity, again derived from the orginal terminal decision problem, which may likewise be used to formulate the design criterion.

Section 8 considers the relationship between the problems of prediction and estimation, showing that each may be reformulated in terms of the other, and derives an "estimative" loss function which may be used to solve the predictive decision problem.

Thus far the various functions considered have all been derived directly from the original terminal decision problem, and we have shown mutual consistency of the variant forms of the design criterion based on the appropriately derived functions. However, it is possible to begin the analysis with an uncertainty, or discrepancy, or dependence function, and use that as the basis of a design crite-rion. We thus wish to know when such a function is *coherent, i.e.* can be obtained from some decision problem. The rest of the paper tackles these characterisation problems. Section 9 describes a regularity condition which will be needed at several points. Then §10 addresses the characterisation of coherent uncertainty functions. It is shown that, under suitable additional technical conditions, coherence is essentially equivalent to concavity. Section 11 characterises coherent discrepancy functions, and relates them to a similar definition introduced by Aitchison. In particular, this refutes Aitchison's conjecture that the only discrepancy function satisfying his definition is the Kullback-Liebler distance. Section 12 gives some properties of coherent dependence funtions, although it stops short of a complete characterisation. Section 13 collects together the various relationships established between the different functions discussed. Once again, we shall have mutual consistency of the different design criteria, when they are based on appropriately related functions. Finally, §14 offers some conclusions.

Throughout the paper we illustrate the theory with several running examples of independent interest. In particular, Example 2 establishes connexions between the log score, entropy, Kullback-Liebler distance, and cross-entropy; while Example 8 brings the familiar optimal design criterion of $D$-optimality under the umbrella of Bayesian decision theory.

## 2 Predictive experimental design

The relationships and results to be developed in this paper, although useful for a wide variety of purposes, are particularly well illustrated by their application to the general problem of *predictive experimental design*, to which we shall return repeatedly. A more specific study of the connexions between our decision-theoretic approach and the "classical" theory of optimal design for normal linear models (amongst others) is undertaken by Dawid and Sebastiani (1996).

Suppose, then, that we shall ultimately wish to predict a quantity $Z \in \mathcal{Z}$. In decision-theoretic terms, this may be expressed as a requirement to select some action $a \in \mathcal{A}$, where $\mathcal{A}$ may be a space of point predictions or an entirely general *action space*. The consequence of this choice, if the eventual value of $Z$ is $z$, is measured by the ensuing *loss* $L(z, a)$. This set-up will be called the *terminal decision problem.*

We suppose given a statistical model for $Z$, such that $Z$ has distribution $P_\theta^Z$, conditional on the value $\theta$ of a parameter $\Theta$. (In general, we shall annotate $P$ with a superscript naming the variable whose distribution is being considered, and with subscripts giving the values of any variables being conditioned on; consequently a random variable as a subscript will imply that we are dealing with a random distribution.) There is also a class $\Xi$ of possible experiments we can perform to obtain information about $\Theta$—and thus, indirectly, about $Z$. If $\xi \in \Xi$ is performed, then observation will be made of a quantity $Y \equiv Y_\xi$, supposed to have distribution

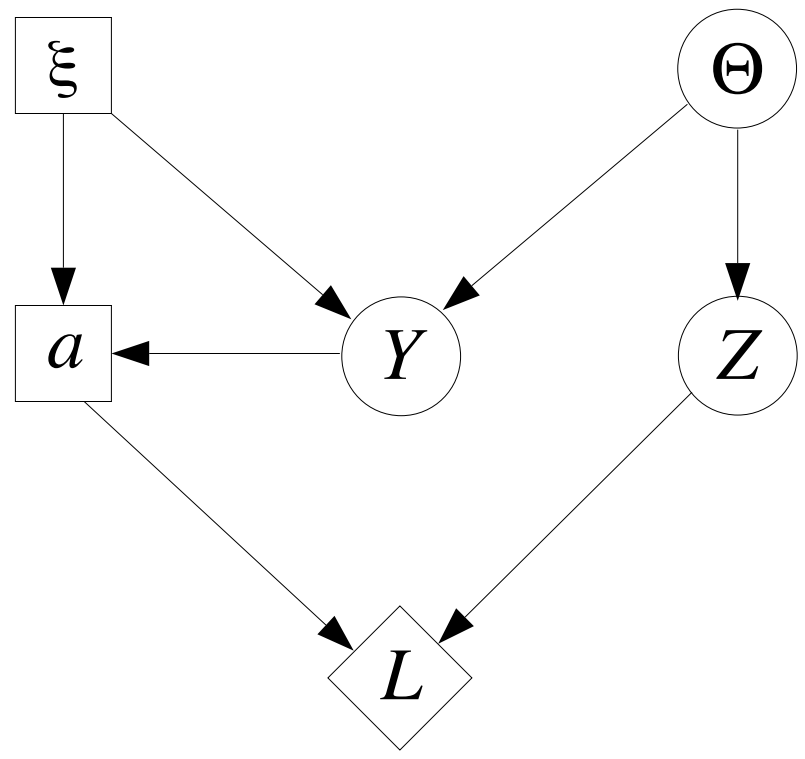

Figure 1: Influence diagram for predictive design

$$
\begin{array}{ccccccccccccc}
 & \xi & & y & & a & & \theta & & z & & & \\
\square & \rightarrow & \bigcirc & \rightarrow & \square & \rightarrow & \bigcirc & \rightarrow & \bigcirc & \rightarrow & \bullet & & L(z,a) \\
\nu_0 & & \nu_1 & & \nu_2 & & \nu_3 & & \nu_4 & & \nu_5 & &
\end{array}
$$

Figure 2: Path through decision tree for predictive design

$P^Y_{\theta,\xi}$ over a sample space $\mathcal{Y}_\xi$ when $\Theta = \theta$. After observing $Y = y$, we then face the terminal decision problem: an action $a \in \mathcal{A}$ is chosen, then $Z$ is observed, and finally the loss suffered if $Z = z$ is $L(z, a)$.

It is supposed that $Z \perp\!\!\!\perp (\xi, Y, a) \mid \Theta$; that is, for each experiment $\xi$, given $\Theta = \theta$, $Z$ is independent of $Y$, with distribution unaffected by the choice of $\xi \in \Xi$ and $a \in \mathcal{A}$. To complete the specification, we assign to $\Theta$ a fixed prior distribution $\Pi$, independent of the choice of $\xi \in \Xi$. The assumed structure of the problem is thus as displayed in the influence diagram of Figure 1.

Note particularly that the above predictive formulation includes, as a special case, the simpler problem of estimation: we may just take $Z \equiv \Theta$, and thus $P^Z_\theta = \delta_\theta$, the one-point distribution at $\theta$. As the additional generality of the prediction problem is easily handled, and of independent interest, we have chosen to develop that here.

Our objective is, first to choose the experiment $\xi$ from $\Xi$, and then, after observing $Y$, to choose the optimal action $a$ from $\mathcal{A}$, in such a way as to minimise the overall expected loss. The emphasis here will be on the choice of experiment $\xi$.

A typical path through the decision tree for this problem looks like Figure 2. We can solve this using the standard extensive analysis method of "backwards induction" (see, *e.g.*, Raiffa and Schlaifer, 1961, §1.2). This proceeds by starting at terminal nodes, $\nu_5$, which are evaluated by their attached loss, and proceeding backwards to evaluate each earlier node: at a random (round) node we do this by taking an expectation over the values of nodes reachable at the subsequent step, using the distribution for that step conditional on the past at that point; and at a decision (square) node we minimise over the immediately reachable values. The decisions yielding the minima at decision nodes then determine the overall optimal strategy. This strategy will comprise (i) specification of the optimal experiment to perform at the initial node $\nu_0$; and (ii) a rule determining, for each node of the form $\nu_2$, the optimal act $a$ to take in the light of observation of data $y$ from experiment $\xi$.

In the sequel we shall consider this procedure in more detail, and, in particular, show how the ensuing criterion for choosing $\xi$ in (i) above may be re-expressed in a variety of interesting ways, using general concepts and measures of properties of distributions which can be associated with the terminal decision problem. These include:

- scoring rules, which encourage honest assessment of probability distributions;
- uncertainty functions, measuring the lack of precision in a distribution;
- discrepancy functions, measuring the "distance" between two distributions; and
- dependence functions, measuring the dependence between two random variable.

Moreover, although we shall here concentrate on their use in experimental design, all these concepts have their own intuitive meaning, and a variety of other applications.

# 3 Terminal decision problem

We start by considering a general terminal decision problem, where the consequence depends on the decision-maker's choice of an act $a$ from the action space $\mathcal{A}$, and on Nature's choice of a value $z$ for an uncertain quantity $Z$ with values in $\mathcal{Z}$. The axiomatic approach to decision-making (see, *e.g.*, DeGroot, 1970; Fishburn, 1982) deduces, from desirable properties of rational attitudes to choices under uncertainty, the existence of a loss function (or negative utility function) $L(z, a)$, and of a distribution $P$ for $Z$, such that it is appropriate to evaluate any act $a$ in terms of its associated *expected loss*

$$L(P, a) := E_P \{L(Z, a)\}, \tag{3.1}$$

so long as this exists.

We are here concerned with problems where the loss function $L$, representing preferences among consequences, is fixed; but the distribution $P$, representing uncertainty about $Z$, is allowed to vary, to take account of changes in relevant information. In a statistical problem, $P$ could represent prior uncertainty about $Z$, in the absence of data; or posterior uncertainty, taking into account any data that may have been gathered. In general, we shall initially restrict $P$ to some convex family $\mathcal{P}_0$ of distributions over $\mathcal{Z}$, which we term the *pre-domain* of the problem. Here and throughout the paper we understand convexity to mean that, if $\tilde{P}$ is a random element of $\mathcal{P}_0$, then $E(\tilde{P}) \in \mathcal{P}_0$. We shall require that, for all $P \in \mathcal{P}_0$, $a \in \mathcal{A}$,

$$L(P, a) \text{ exists, as a finite value or } \pm\infty. \tag{3.2}$$

For $P \in \mathcal{P}_0$, define $H : \mathcal{P}_0 \to [-\infty, \infty]$ by:

$$H(P) := \inf_{a \in \mathcal{A}} L(P, a). \tag{3.3}$$

Then $H(P)$ is the *Bayes loss* under $P$. It is readily seen that $H$ is concave over $\mathcal{P}_0$.

We shall need further to restrict attention to a convex subset $\mathcal{P} \subseteq \mathcal{P}_0$, the *domain* of the problem (or of $L$), such that, for any $P \in \mathcal{P}$, $H(P)$ is finite. We shall also suppose throughout that the infimum in (3.3) is actually achieved, at a (not necessarily unique) *Bayes act*, $a(P) \in \mathcal{A}$ (however, we may allow $L(P, a)$ to be $+\infty$ for some $a \in \mathcal{A}$). Thus $a(P)$ is an optimal act when uncertainty is expressed by $P \in \mathcal{P}$.

When $L$ is bounded, we can take as domain the set of all distributions over $\mathcal{Z}$. More generally, suppose that, for all $P \in \mathcal{P}_0$, $H(P) < \infty$. Then we can take as domain $\mathcal{P} := \{P \in \mathcal{P}_0 : H(P) > -\infty\}$, and this is convex, by concavity of $H$. This will apply in particular if, for all $P \in \mathcal{P}_0$, $H(P) \leq 0$, in which case we shall call the problem *regular* with respect to $\mathcal{P}_0$.

Note that, for $\tilde{P}$ a random element of $\mathcal{P}$, $E\{H(\tilde{P})\}$ exists and is finite, by concavity and finiteness of $H$ on $\mathcal{P}$.

If, for $P, Q \in \mathcal{P}, P \neq Q$, Bayes acts under $P$ and $Q$ are necessarily different, we shall call $L$ *strict*. Evidently, this property can only hold when the cardinality of the action space $\mathcal{A}$ is at least as large as that of $\mathcal{P}$.

## 3.1 Equivalence

Suppose that a new loss function $L^*$ is defined by:

$$L^*(z, a) \equiv cL(z, a) + k(z), \tag{3.4}$$

where $c > 0$ and $k$ is an arbitrary function. Then, defining

$$k(P) := E_P\{k(Z)\}, \tag{3.5}$$

so long as $k(P)$ is finite we have

$$L^*(P, a) \equiv cL(P, a) + k(P). \tag{3.6}$$

Under (3.4), $L^*$ and $L$ yield the same comparisons among acts: $L^*(z,a) - L^*(z,a') \equiv c\{L(z,a) - L(z,a')\}$. When $k(P)$ is finite, $L^*(P,a) - L^*(P,a') \equiv c\{L(P,a) - L(P,a')\}$. In particular, the optimal act $a(P)$ is then the same, whether $L^*$ or $L$ is the loss function.

We shall call two loss functions *equivalent* if they are related as in (3.4), and *strongly equivalent* if moreover $c = 1$. If $k$ is bounded, we can take $L$ and $L^*$ to have the same domain, but otherwise the domain of $L^*$ may need to differ from that of $L$, so as to ensure the finiteness of $H^*(P) := \inf_{a \in \mathcal{A}} L^*(P,a)$.

In particular, fix $a_0 \in \mathcal{A}$, and take $c = 1$, $k(z) \equiv -L(z, a_0)$. Then clearly $H^*(P) \leq 0$. Thus, by an equivalence transformation, we may always take the uncertainty function to be non-positive. The transformed problem is thus regular with respect to any convex pre-domain $\mathcal{P}_0^*$ for which (3.2) holds for $L^*$, and consequently we may take its domain as $\mathcal{P}^* := \{P \in \mathcal{P}_0^* : H^*(P) > -\infty\}$.

## 3.2 Uncertainty functions

We have noted that $H$ is a concave function of $P$ (DeGroot, 1970, §8.4). Various authors (*e.g.* DeGroot, 1962; Rao, 1982) have used general concave functions of $P$ as measures of the uncertainty or diversity in $P$. Any function of $P$ which can arise as the Bayes loss in some terminal decision problem will be termed a *coherent uncertainty function.* In §10 we shall address the characterisation of such functions.

When $L^*$ is equivalent to $L$, as in (3.4), we likewise shall term the corresponding uncertainty functions $H$ and $H^*$ *equivalent* (*strongly equivalent* if $c = 1$). Then

$$H^*(P) \equiv cH(P) + k(P) \tag{3.7}$$

on the intersection of their domains.

For the case that $L$ is bounded, with $\mathcal{P}$ containing all distributions over $\mathcal{Z}$, an uncertainty function $H$ may be called *canonical* if $H(P) \geq 0$, with equality when $P$ is a one-point distribution. If $H$ is not canonical, then its strongly equivalent version $H^*$, given by

$$H^*(P) \equiv H(P) - E_P\{H(\delta_Z)\}, \tag{3.8}$$

where $\delta_z$ denotes the one-point distribution at $z \in \mathcal{Z}$, will be so.

**Example 1** Let $t : \mathcal{Z} \to \mathbb{R}^1$, $\mathcal{A} = \mathbb{R}^1$ and take $L(z,a) \equiv \{t(z) - a\}^2$, corresponding to the problem of estimating $T := t(Z)$, with quadratic loss. This was used as a loss function for predictive experimental design by Lindley (1968), Brooks (1972). We can take as pre-domain $\mathcal{P}_0$ the set of all distributions $P$ for $Z$ such that $E_P(T)$ exists in $[-\infty, \infty]$. Then $L(P,a) \equiv \mathrm{var}_P(T) + \{E_P(T) - a\}^2$, and we obtain $a(P) \equiv E_P(T)$, $H(P) \equiv \mathrm{var}_P(T)$, defined over domain $\mathcal{P} := \{P : \mathrm{var}_P(T) < \infty\}$. This $H$ is canonical. Unless $\#\{t(\mathcal{Z})\} \leq 2$, $L$ is not strict.

Alternatively, we can replace $L$ by the strongly equivalent loss function $L^*(z,a) := L(z,a) - L(z,0) \equiv -2at(z) + a^2$. Then $L^*(P,a) \equiv -2aE_P(T) + a^2$. Once again, the Bayes act is $a(P) \equiv E_P(T)$, but the domain may now extended to $\mathcal{P}^* := \{P : E_P(T) \text{ is finite}\}$. The associated (regular) uncertainty function now becomes $H^*(P) \equiv -E_P(T)^2$.

For $P \notin \mathcal{P}^*$, it becomes impossible to compare different acts, and the Bayes act is undefined. □

**Example 2** Let $\mathcal{M}$ be the set of distributions over $\mathcal{Z}$ which are mutually absolutely continuous with respect to a probability measure $\mu$. Take both action space $\mathcal{A}$ and pre-domain $\mathcal{P}_0$ to be $\mathcal{M}$, and define $L(z,Q) \equiv -\log q(z)$, where $q$ is some fixed version of the density of $Q$ with respect to $\mu$. This is the negative log-likelihood, also called the logarithmic score. The expected loss is

$$L(P,Q) \equiv -\int p(z) \log q(z) d\mu, \tag{3.9}$$

where $p = dP/d\mu$. Then $K(P,Q) := L(P,Q) - L(P,P)$ is just the *Kullback-Liebler discrepancy* between $P$ and $Q$:

$$K(P,Q) \equiv \int p(z) \log\{p(z)/q(z)\} d\mu. \tag{3.10}$$

Formulae (3.9) and (3.10) above may be applied more generally: so long as $P \ll Q \ll \mu$, $L(P,Q)$ and $K(P,Q)$ are well-defined, in the sense that they will not depend on the versions of the density $p \equiv dP/d\mu$ and $q \equiv dQ/d\mu$ employed. However the value of $L$ will depend on the choice of base measure $\mu$. As for $K(P,Q)$, its value is, further, independent of $\mu$, so that it is well-defined for any pair $P$, $Q$ with $P \ll Q$ (and we could take, for instance, $\mu = Q$).

For $P \ll Q$ we have, by a well-known inequality,

$$K(P,Q) \geq 0, \tag{3.11}$$

with $K(P,Q) = 0$ holding if and only if $Q = P$.

Since, for $P \in \mathcal{M}$,

$$L(P,P) \equiv -K(P,\mu) \leq 0, \tag{3.12}$$

the problem is regular, and we thus take as domain $\mathcal{P} := \{P \in \mathcal{M} : K(P,\mu) < \infty\}$, which is convex. From (3.11), for fixed $P \in \mathcal{P}$, $L(P,Q) \geq L(P,P)$, all $Q \in \mathcal{P}$, so that $Q(P) = P$ gives the unique Bayes act in this problem. The associated uncertainty function is then

$$H(P) \equiv -\int p(z) \log p(z) d\mu, \tag{3.13}$$

the *entropy* of the distribution $P$ (with respect to the distribution $\mu$).

Changing $\mu$ to an equivalent probability measure $\nu$ results in a strongly equivalent loss function, but may require a change in the domain. In particular, suppose we use as $\nu$ some fixed distribution $R \in \mathcal{M}$, having density $r$ with respect to $\mu$. This yields new loss function

$$L^*(z, Q) := -\log\{q(z)/r(z)\} \equiv L(z, Q) + \log r(z).$$

Then

$$\begin{aligned} L^*(P, Q) &\equiv -\int p(z) \log\{q(z)/r(z)\} d\mu, \\ L^*(P, P) &\equiv -K(P, R). \end{aligned}$$

Once again, (3.11) yields $Q(P) = P$ as the Bayes act. The new domain is $\mathcal{P}^* := \{P : K(P, R) < \infty\}$, which is convex, and the corresponding regular uncertainty function over $\mathcal{P}^*$ is now $H^*(P) \equiv -K(P, R)$.

Note that, for any $P \in \mathcal{M}$, we can always find a distribution $R \in \mathcal{M}$, to use as base measure $\nu$, such that $K(P, R) < \infty$—for example, one possibility is to take $R = P$, for which $K(P, R) = 0$. Thus, for $Z \sim P$, we can always, by some choice of base measure, and thus of domain, compare the act $P$ with any other act $Q$, so validating the conclusion that, for any $P \in \mathcal{M}$, $Q(P) = P$ is the Bayes act. However, it is not in general possible to find a base measure which works in this way for every $P \in \mathcal{M}$ simultaneously, since (for continuous $\mu$, at least) it may be shown that, for any $R \in \mathcal{M}$, there exists $Q \in \mathcal{M}$ with $K(Q, R) = \infty$. Thus we can not take the whole of $\mathcal{M}$ as domain: there is no version of the uncertainty function which can be applied to the whole of $\mathcal{M}$, since any such would take the value $-\infty$ somewhere.

The logarithmic score and entropy function are often defined in terms of an arbitrary $\sigma$-finite measure $\mu$, rather than a probability measure. Much of the above analysis still applies, except that we can no longer assert $K(P, \mu) \geq 0$, and so the problem may not be regular. In particular, it is possible that $\{P \in \mathcal{M} : K(P, \mu) < \infty\}$ may not provide a convex domain. However, we can always, by an equivalence transformation, change the base measure to a probability distribution, thus recovering the set-up considered above. □

# 4 Predictive experimental design based on uncertainty functions

We now return to the experimental design framework of §2. We assume that, for each $\xi \in \Xi$, almost surely (as a function of $Y \equiv Y_\xi$) the predictive distribution $P^Z_{\xi,Y} \in \mathcal{P}$. Then $H(P^Z_{\xi,Y})$ exists almost surely, and has finite expectation over the distribution of $Y$ in $\xi$.

Applying the method of backwards induction to the decision tree of Figure 2 involves the following steps:

($i$) At $\nu_4$ we evaluate the expected loss as

$$E\{L(Z,a)|\xi,y,a,\theta\} \equiv E\{L(Z,a)|\theta\} \equiv L(P^Z_\theta, a),$$

where we have utilised the conditional independence assumption $Z \perp\!\!\!\perp (\xi, Y, a) \mid \Theta$.

($ii$) At $\nu_3$ we then calculate $E\{L(P^Z_\Theta, a)|\xi,y,a\}$. (Recall that $P^Z_\Theta$ is a function of the random conditioning variable $\Theta$, and hence itself random.) This reduces to $L(P^Z_{\xi,y}, a)$, where $P^Z_{\xi,y}$ denotes the predictive distribution for $Z$ after observing $Y = y$ in $\xi$. (This could alternatively have been obtained directly from a simplified tree in which node $\nu_4$, and the arrow into it, was omitted. However, it is helpful to retain the more detailed structure of Figure 2.)

($iii$) At $\nu_2$ we select $a$ to minimise $L(P^Z_{\xi,y}, a)$, *viz.* we take $a = a(P^Z_{\xi,y})$.

At this point we note that the value to be attached to $\nu_2$, the minimised expected loss, is $H(P^Z_{\xi,y})$, where $H$ is the uncertainty function associated with $L$.

($iv$) At $\nu_1$ the expected loss is thus $E\{H(P^Z_{\xi,Y})|\xi\}$, the expectation being over the marginal distribution of $Y$ in $\xi$.

($v$) Finally, at $\nu_0$, the experiment $\xi$ is chosen to minimise

$$E\{H(P^Z_{\xi,Y})|\xi\}. \tag{4.1}$$

Thus we see that the criterion for choosing $\xi$ can be expressed as:

*Minimise the expected predictive uncertainty about $Z$,*
*posterior to observing the outcome $Y$ of $\xi$.*

The expectation here is with respect to the marginal distribution for $Y$ in experiment $\xi$, constructed by combining the conditional distributions $P^Y_{\xi,\theta}$ specified by $\xi$ with the fixed prior distribution $\Pi$ for $\Theta$. When the conditioning on $\xi$ in expressions such as (4.1) is clear from the context, we shall drop it from the notation.

Note that, if $H^*$ is equivalent to $H$, satisfying (3.7) with $k(P^Z)$ finite, then

$$E\left\{H^*\left(P^Z_{\xi,Y}\right)\middle|\,\xi\right\} = cE\left\{H\left(P^Z_{\xi,Y}\right)\middle|\,\xi\right\} + k\left(P^Z\right),$$

so that equivalent problems will yield the same choice of experiment at $\nu_0$ (as well as the same optimal act at $\nu_2$).

**Cost function** As a slight extension of the above analysis, we might also have an additive cost $c(\xi, y)$ associated with deciding to perform the experiment $\xi$ and then observing $Y = y$ (often $c(\xi, y)$ will not depend on $y$). It is of course important that $c$ is measured on a scale commensurate with that of $L$, and hence of $H$, so that it is meaningful to form the overall terminal loss $c(\xi, y) + L(z, a)$. It now becomes necessary to add to the criterion (4.1) the quantity $c(\xi) := E\{c(\xi, Y)\}$ (assumed to exist). We shall not take any further explicit account of this simple extension.

**Example 1** In Example 1, the design criterion becomes:

$$\textit{Minimise } E\{\mathrm{var}(T|Y)\}.$$

Alternatively, using the strongly equivalent loss function $L^*$ yields:

$$\textit{Maximise } E\{E(T|Y)^2\}.$$

We note that the sum of these two variant criteria is $E(T^2)$, which does not depend on $\xi$ (as may be confirmed by inspection of Fig. 1, $T$ being a function of $Z$), thus confirming their equivalence when $E(T^2)$ is finite. □

**Example 2** In Example 2, the rule for choosing the design is:

*Minimise the expected entropy of the*
*conditional distribution of $Z$ given $Y$.*

Here we could use the entropy (3.13) with respect to the base measure $\mu$, or, equivalently, the entropy with respect to a distribution $R \in \mathcal{M}$, *viz.* $-K(P, R)$. □

# 5 Scoring rules

We now return to the general structure of §3. Suppose You have to quote a distribution $Q$ for $Z$, from some convex family $\mathcal{P}$ of such distributions. If $Z = z$, You will suffer a loss $S(z, Q)$. This is a particular form of terminal decision problem in which the action space $\mathcal{A}$ is $\mathcal{P}$. Writing $S(P, Q) := E_P\{S(Z, Q)\}$, such a "scoring rule" encourages honesty if

$$S(P, P) \leq S(P, Q), \quad \text{all} \quad P, Q \in \mathcal{P}. \tag{5.1}$$

In this case, if we restrict attention to distributions $P \in \mathcal{P}$, the optimal act when $Z \sim P$ is just to quote the true distribution $P$. The corresponding uncertainty function is thus

$$H(P) \equiv S(P, P). \tag{5.2}$$

A scoring rule $S(z, Q)$ satisfying (5.1), with $S(P, P)$ finite, is termed *proper* (over $\mathcal{P}$); it is *strictly proper* if inequality holds in (5.1) whenever $Q \neq P$. Proper scoring rules are an important tool of Bayesian decision theory and inference: see Good (1952), Savage (1971), de Finetti (1975).

Now suppose we start with a general terminal decision problem, as in §3, with loss function $L$ and domain $\mathcal{P}$. Define, for $Q \in \mathcal{P}$,

$$S(z, Q) \equiv L\{z, a(Q)\}. \tag{5.3}$$

(If the Bayes act $a(Q)$ is not unique, some particular choice must be made). Then, for $P, Q \in \mathcal{P}$, $S(P, P) = L\{P, a(P)\} \leq L\{P, a(Q)\} = S(P, Q)$, by definition of $a(P)$. Hence (Dawid, 1986) (5.3) defines a proper scoring rule over $\mathcal{P}$; it will be strictly proper if $L$ is strict. Moreover, the uncertainty function $H(P) \equiv L\{P, a(P)\} \equiv S(P, P)$ will be the same, whether based on $L$ or on $S$ as loss function. In this way, whatever the original decision problem, it may always be reinterpreted as involving the assessment of a probability distribution for $Z$.

**Equivalence** If $L^*$ is equivalent to $L$ as in (3.4), with associated scoring rule $S^*$, then

$$S^*(z, Q) \equiv cS(z, Q) + k(z).$$

In this case we may call $S^*$ and $S$ *equivalent* (*strongly equivalent* if $c = 1$). Then

$$S^*(P, Q) \equiv cS(P, Q) + k(P),$$

when $k(P) \equiv E_P\{k(Z)\}$ exists.

**Example 1** We have

$$S(z, Q) \equiv \{t(z) - E_Q(T)\}^2, \tag{5.4}$$

$$S(P, Q) \equiv \mathrm{var}_P(T) + \{E_P(T) - E_Q(T)\}^2. \tag{5.5}$$

In general, this scoring rule is not strict, since $S(P, Q) = S(P, P)$ whenever $E_Q(T) = E_P(T)$. □

**Example 2** We note that the original loss function for Example 2 is already a strictly proper scoring rule over $\mathcal{P}$. Whenever our original decision problem has $\mathcal{A} = \mathcal{P}$ and $L$ is strictly proper, (5.3) simply recovers $S \equiv L$.

In this example we have

$$S(P, Q) \equiv -\int p(z) \log q(z) d\mu \qquad (P, Q \in \mathcal{P}); \tag{5.6}$$

or, using the strongly equivalent form $L^*$, with base measure $R \in \mathcal{M}$,

$$S^*(P, Q) \equiv -\int p(z) \log\{q(z)/r(z)\} d\mu \qquad (P, Q \in \mathcal{P}^*). \tag{5.7}$$

□

$$\begin{array}{cccccccccccc} & \xi & & y & & Q & & \theta & & z & & \\ \square & \rightarrow & \bigcirc & \rightarrow & \square & \rightarrow & \bigcirc & \rightarrow & \bigcirc & \rightarrow & \bullet & S(z,Q) \\ \nu_0 & & \nu_1 & & \nu_2 & & \nu_3 & & \nu_4 & & \nu_5 & \end{array}$$

Figure 3: Decision tree for probability prediction

### 5.1 Predictive design using scoring rules

Suppose that, in the predictive design problem of §2, we replace $\mathcal{A}$ by $\mathcal{P}$, $a \in \mathcal{A}$ by $Q \in \mathcal{P}$, and $L(z,a)$ by the proper scoring rule $S(z,Q)$ associated with $L$. This does not affect the uncertainty function $H$. The decision tree is now as in Figure 3. Applying the analysis of §4 to this new problem, we have, for the initial steps:

$(i)'$ At $\nu_4$, find $S(P^Z_\theta, Q)$.

$(ii)'$ At $\nu_3$, find $S(P^Z_{\xi,y}, Q)$.

$(iii)'$ At $\nu_2$, choose $Q = P^Z_{\xi,y}$.

At this point we once again obtain $H(P^Z_{\xi,y})$ as the value to be attached to $\nu_2$. Consequently, the remaining steps, $(iv)'$ and $(v)'$, will agree with $(iv)$ and $(v)$ of §4, so that the optimal design is the same as before. It follows that we may, without loss of generality, always assume that our terminal decision problem is one of probability prediction, using a proper scoring rule $S$ for $Z$. The optimal design is again given by minimising (4.1), for the coherent uncertainty function $H(P) \equiv S(P,P)$.

**Example 1** In particular, in Example 1, instead of using the loss function $L(z,a) \equiv \{t(z) - a\}^2$ ($a \in \mathbb{R}^1$), it is equivalent to use, as loss function, the proper scoring rule $S(z,Q) \equiv \{t(z) - E_Q(T)\}^2$. In either case, the design criterion is:

$$\textit{Minimise } E\{\text{var}(T \mid Y)\},$$

the expectation being over the distribution of $Y$ in $\xi$. □

### 5.2 Further proper scoring rules

A wide variety of proper scoring rules has been proposed. Here we present some examples of particular interest, noting their associated design criteria.

We first note the following. Let $\{S_\lambda : \lambda \in \Lambda\}$ be a family of proper scoring rules, and let $W$ be a measure on $\Lambda$. Define

$$S(z,Q) := \int S_\lambda(z,Q)dW(\lambda), \tag{5.8}$$

assuming that the integral exists. Then $S$ is also proper, and, under conditions allowing the interchange of integration and expectation, $S(P,Q) \equiv \int S_\lambda(P,Q)dW(\lambda)$; while the associated uncertainty functions similarly satisfy $H(P) \equiv \int H_\lambda(P)dW(\lambda)$.

**Example 3** Let $t$ be a complex function on $\mathcal{Z}$, $T \equiv t(Z)$, and define

$$S(z,Q) \equiv |t(z) - E_Q(T)|^2. \tag{5.9}$$

Separating real and imaginary parts, this is seen to be a combination of two cases of (5.4), hence proper. We obtain

$$\begin{aligned} S(P,Q) &\equiv \text{var}_P(T) + |E_P(T) - E_Q(T)|^2, & (5.10)\\ H(P) &\equiv \text{var}_P(T), & (5.11) \end{aligned}$$

where $\text{var}_P(T)$ is now interpreted as $E_P(|T|^2) - |E_P(T)|^2$. As defined, the appropriate domain is the set of distributions such that $\text{var}_P(T)$ exists; again, by an equivalence transformation, it is enough that $E_P(T)$ should.

The design criterion generated by this scoring rule is, as in Example 1:

$$\textit{Minimise } E\{\text{var}(T \mid Y)\},$$

but now using the above extended definition of "var".

In particular, with $\mathcal{Z} = \mathbb{R}^m$, take, for fixed $u \in \mathbb{R}^m$, $t(z) \equiv e^{iu^{\mathrm{T}}z}$. We obtain

$$S(z,Q) \equiv \left|e^{iu^{\mathrm{T}}z} - \phi_Q(u)\right|^2, \tag{5.12}$$

where $\phi_Q$ is the characteristic function of $Q$, and

$$\begin{aligned} S(P,Q) &\equiv 1 - |\phi_P(u)|^2 + |\phi_P(u) - \phi_Q(u)|^2, & (5.13)\\ H(P) &\equiv 1 - |\phi_P(u)|^2. & (5.14) \end{aligned}$$

Since $S$ is bounded, there are no difficulties about existence of expectations, and we may take as domain the space of all distributions for $Z$.

The design criterion in this case can be expressed as:

$$\textit{Maximise } E\left|\phi_Y^Z(u)\right|^2,$$

where $\phi_y^Z$ denotes the characteristic function of the conditional (predictive) distribution, $P_y^Z$, of $Z$ given observation of $y$ (in experiment $\xi$); and the expectation is with respect to the distribution of $Y$ in $\xi$.

Note that $S$ in (5.12) is invariant under the translations $z \to z+a$, $Q \to Q*\delta_a$, $a \in \mathbb{R}^m$, where $\delta_a$ is the point mass at $a$ and $*$ denotes convolution. Likewise, (5.13) is invariant under $P \to P * \delta_a$, $Q \to Q * \delta_a$. □

**Example 4** With $\mathcal{Z} = \mathbb{R}^m$, it now follows from (5.8) and (5.12) that a proper scoring rule is given by

$$S(z, Q) \equiv \int \left| e^{iu^{\mathrm{T}}z} - \phi_Q(u) \right|^2 dW(u), \tag{5.15}$$

yielding

$$\begin{aligned} S(P, Q) &\equiv \int \left\{ 1 - |\phi_P(u)|^2 \right\} dW(u) \\ &\qquad + \int |\phi_P(u) - \phi_Q(u)|^2 \, dW(u), \end{aligned} \tag{5.16}$$

$$H(P) \equiv \int \left\{ 1 - |\phi_P(u)|^2 \right\} dW(u). \tag{5.17}$$

If the support of $W$ is $\mathbb{R}^m$, then $S$ is strictly proper.

The corresponding design criterion is:

$$\textit{Maximise} \int E \left| \phi_Y^Z(u) \right|^2 dW(u).$$

□

**Example 5** (Eaton, 1982; Eaton *et al.*, 1996). Let $K : \mathcal{Z} \times \mathcal{Z} \to \mathbb{C}$ be a non-negative definite kernel: thus $K(a, b) \equiv \overline{K(b, a)}$, and $\|\alpha\|_K^2 := \int \int K(a, b) d\alpha(a) d\alpha(b) \geq 0$ for all bounded signed measures $\alpha$ on $\mathcal{Z}$. Take $\mathcal{P} = \{P : \|P\|_K^2 < \infty\}$ and, for $Q \in \mathcal{P}$, define

$$S(z, Q) := \|Q\|_K^2 - \int K(z, b) dQ(b) - \int K(a, z) dQ(a), \tag{5.18}$$

strongly equivalent to $\|Q - \delta_z\|_K^2$. For the case that $\mathcal{Z}$ is finite, with $K$ the Kronecker delta function, $\|\alpha\|$ is just the Euclidean length of the mass function (density) of $\alpha$, and $\|Q - \delta_z\|_K^2$ then defines the "quadratic score", or the *Brier score* (see Dawid, 1986) in the binary case.

We find

$$S(P, Q) \equiv \int \int K(a, b) \left\{ dQ(a) dQ(b) - dP(a) dQ(b) - dQ(a) dP(b) \right\}, \tag{5.19}$$

so that

$$H(P) \equiv -\|P\|_K^2 \tag{5.20}$$

and so

$$S(P, Q) - S(P, P) \equiv \|P - Q\|_K^2. \tag{5.21}$$

Thus (5.18) defines a proper scoring rule, which is strictly proper if $K$ is positive definite.

The design criterion associated with (5.18) becomes:

$$Maximise\ E\left\{\left\|P_Y^Z\right\|_K^2\right\}.$$

When $K(a,b) \equiv t(a)\overline{t(b)}$, (5.18) is equivalent to (5.9). In particular, for $K(a,b) \equiv e^{iu^{\mathrm{T}}(a-b)}$, (5.18) is equivalent to (5.12), and thus we recover (5.15) as a special case for $K(a,b) \equiv \int e^{iu^{\mathrm{T}}(a-b)}dW(u)$. □

**Example 6** As a variant of Example 5, let $\mu$ be a base measure on $\mathcal{Z}$, $\mathcal{M}$ the set of signed measures on $\mathcal{M}$ absolutely continuous with respect to $\mu$, and $\mathcal{P} \subset \mathcal{M}$ its subset of probability measures. The norm $\|A\|_\mu$ of $A \in \mathcal{M}$ is now defined by: $\|A\|_\mu^2 = \int a(z)^2 d\mu(z)$, where $a \equiv dA/d\mu$. Now we take

$$S(z,Q) \equiv \|Q\|_\mu^2 - 2q(z),$$

generalizing the "quadratic score" of Example 5.

We obtain

$$\begin{aligned} S(P,Q) &\equiv \|P-Q\|_\mu^2 - \|P\|_\mu^2, \\ H(P) &\equiv -\|P\|_\mu^2, \end{aligned}$$

and we again see that $S$ is a proper scoring rule.

Now the design criterion is:

$$Maximise\ E\left\{\left\|P_Y^Z\right\|_\mu^2\right\}.$$

□

**Example 7** Let $\mathcal{E}$ be a regular full exponential family, having canonical statistic $T \equiv t(Z)$ with values in $\mathbb{R}^k$, and parametrised by the expectation parameter $\tau := E(T)$, taking values in $\mathcal{T}$. The density with respect to some underlying measure $\nu$ on $\mathcal{Z}$ is thus of the form

$$f(z;\tau) \equiv \exp\left\{a(z) - b(\lambda) + \lambda^{\mathrm{T}}t(z)\right\}, \tag{5.22}$$

where $\beta(\lambda) = \tau$, with $\beta := \nabla b$, defines the $(1,1)$ relationship between $\tau$ and the natural parameter $\lambda$. Let $\mathcal{P}$ be the family of distributions $P$ for $Z$ such that $\tau_P := E_P(T)$ exists and lies in $\mathcal{T}$. For $Q \in \mathcal{P}$ define

$$S(z,Q) := -\log f(z;\tau_Q). \tag{5.23}$$

Note that (5.23) coincides with the logarithmic score for $Q \in \mathcal{E}$, but not outside $\mathcal{E}$.

Defining $\lambda_P$ by $\beta(\lambda_P) = \tau_P$, we have, up to strong equivalence:

$$\begin{aligned} S(z,Q) &\equiv b(\lambda_Q) - \lambda_Q^{\mathrm{T}}t(z), \\ S(P,Q) &\equiv b(\lambda_Q) - \lambda_Q^{\mathrm{T}}\tau_P, \end{aligned}$$

so that

$$S(P,Q) - S(P,P) \equiv \{b(\lambda_Q) - b(\lambda_P)\} - (\lambda_Q - \lambda_P)^{\mathrm{T}}\tau_P, \tag{5.24}$$

which, since $\tau_P = (\nabla b)(\lambda_P)$, is non-negative by convexity of $b$. Hence $S$ is proper, with associated uncertainty function

$$H(P) \equiv b(\lambda_P) - \lambda_P^{\mathrm{T}}\tau_P. \tag{5.25}$$

However, $S$ is not strictly proper, since $S(P,Q) = 0$ whenever $E_Q(T) = E_P(T)$.

For the case $\mathcal{E} = \{N(\tau, 1)\}$, $\mathcal{P} = \{P : E_P(T) \text{ exists}\}$, we recover Example 1 (up to equivalence). □

**Example 8** As another instance of Example 7, let $\mathcal{E}$ be the family of $p$-variate normal distributions for $Z$, with arbitrary mean vector $\mu$ and dispersion matrix $\Sigma$. We can take $T = (Z, ZZ^{\mathrm{T}})$, $\tau = (\mu, \mu\mu^{\mathrm{T}} + \Sigma)$. Let $\mu_P, \Sigma_P$ be the mean and dispersion of $Z$ under $P \in \mathcal{P}$, where $\mathcal{P}$ is the set of distributions for which these moments exist. We then obtain the proper scoring rule over $\mathcal{P}$:

$$S(z,Q) \equiv \frac{1}{2}\left\{\log\det\Sigma_Q + (z-\mu_Q)^{\mathrm{T}}\Sigma_Q^{-1}(z-\mu_Q)\right\}. \tag{5.26}$$

This yields

$$\begin{aligned} S(P,Q) &\equiv \frac{1}{2}\left\{\log\det\Sigma_Q + (\mu_P-\mu_Q)^{\mathrm{T}}\Sigma_Q^{-1}(\mu_P-\mu_Q)\right. \\ &\qquad \left. + \operatorname{tr}(\Sigma_Q^{-1}\Sigma_P)\right\}, &(5.27) \\ \text{and } H(P) &\equiv \frac{1}{2}\left\{\log\det\Sigma_P + p\right\}. &(5.28) \end{aligned}$$

To verify propriety of $S$ directly, note that

$$\begin{aligned} 2\{S(P,Q) - S(P,P)\} &\geq \operatorname{tr}\Sigma_Q^{-1}\Sigma_P - \log\det\Sigma_Q^{-1}\Sigma_P - p \\ &= \sum_i (\lambda_i - \log\lambda_i - 1) \\ &\geq 0, \end{aligned}$$

where the $(\lambda_i)$ are the eigenvalues of $\Sigma_Q^{-1}\Sigma_P$.

We thus obtain the design criterion:

$$\textit{Minimise } E\left\{\log\det\Sigma_Y^Z\right\},$$

where $\Sigma_y^Z := \Sigma_{P_y^Z}$ denotes the predictive dispersion matrix of $Z$ given $Y = y$ (in $\xi$); and the expectation is over the marginal distribution of $Y$ in $\xi$. For the case of an estimation problem, where $Z$ is identified with the parameter $\Theta$, $\Sigma_y^Z$ is just the posterior dispersion matrix. The resulting design criterion is closely related to the criterion of $D$-optimality used in the general theory of optimal experimental design (Dawid and Sebastiani, 1996). □

# 6 Discrepancy functions

Let $P, Q \in \mathcal{P}$, a suitable family of distributions. We can attempt to measure the discrepancy between these distributions by means of some function $D(P,Q)$, not necessarily symmetric, such that $D(P,Q) \geq 0$, and $D(P,P) \equiv 0$. Such a discrepancy $D$ is termed *strict* if $D(P,Q) > 0$ for $P \neq Q$. In this Section we show how we can naturally associate a discrepancy function with any terminal decision problem.

Indeed, from (5.1), we see that, given any proper scoring rule $S$ on domain $\mathcal{P}$, with associated uncertainty function $H(P)$, we may define a discrepancy on $\mathcal{P}$ by:

$$D(P,Q) \equiv S(P,Q) - S(P,P) \tag{6.1}$$

$$\equiv S(P,Q) - H(P). \tag{6.2}$$

Then $D$ will be strict if $S$ is. In particular, this construction may be initiated with an arbitrary loss function $L$, *via* (5.3), yielding

$$D(P,Q) \equiv L\{P, a(Q)\} - L\{P, a(P)\}. \tag{6.3}$$

A discrepancy function which can arise in this way from a decision problem will be called *coherent*. We shall consider properties and characterisations of coherent discrepancy functions in §11.

It is easily seen that any two strongly equivalent decision problems will yield the identical discrepancy function, although possibly over different domains. We could thus extend $D$ to be defined over the union of the domains of the various strongly equivalent versions of $L$, or of $S$. More generally, if (3.4) holds, we have $D^*(P,Q) \equiv cD(P,Q)$, in which case we may call $D$ and $D^*$ *equivalent*.

## 6.1 Predictive design with discrepancy functions

Let now $D$ be a coherent discrepancy, arising from some loss function $L$, and let $S$ and $H$ be the associated proper scoring rule and uncertainty function. Step $(i)'$ of §5.1 can then be written as:

$(i)''$ At $\nu_4$, find $H(P_\theta^Z) + D(P_\theta^Z, Q)$,

whence, applying the regular backwards induction routine, the next step is equivalent to:

$(ii)''$ At $\nu_3$, find $E\{H(P_\Theta^Z)|\xi, y\} + E\{D(P_\Theta^Z, Q)|\xi, y\}$.

Then:

$(iii)''$ At $\nu_2$ choose, in accordance with $(iii)'$ of §5.1, $Q = P^Z_{\xi,y}$, thus obtaining

$$E\left\{H\left(P^Z_\Theta\right)\middle|\,\xi, y\right\} + E\left\{D\left(P^Z_\Theta, P^Z_{\xi,y}\right)\middle|\,\xi, y\right\}.$$

$(iv)''$ At $\nu_1$, find $E\{H(P^Z_\Theta)\} + E\{D(P^Z_\Theta, P^Z_{\xi,Y})|\xi\}$.

Note that the first term in $(iv)''$ does not depend on $\xi$. So the final step is:

$(v)''$ At $\nu_0$ choose $\xi$ to minimise

$$E\{D(P^Z_\Theta, P^Z_{\xi,Y})|\xi\}. \tag{6.4}$$

That is, the design criterion may be re-expressed as:

*Minimise the expected discrepancy between the sampling distribution of $Z$ and the predictive distribution of $Z$.*

In this criterion, the sampling distribution depends on $\Theta$, and the predictive distribution on $Y$, so that the expectation is over the joint distribution of $(\Theta, Y)$ in $\xi$.

## 6.2 Examples

**Example 1** Here we have

$$D(P, Q) \equiv \{E_P(T) - E_Q(T)\}^2. \tag{6.5}$$

In general this is not strict. In this case the design criterion can thus be re-expressed as:

*Minimise the expected squared difference between the sampling expectation of $T$, and the predictive expectation of $T$.*

□

**Example 2** We have

$$D(P, Q) \equiv K(P, Q), \tag{6.6}$$

the Kullback-Liebler distance $K(P, Q)$ (3.10) between $P$ and $Q$. The design criterion thus becomes:

*Minimise the expected Kullback-Liebler distance of the predictive distribution from the sampling distribution.*

□

**Example 3** The discrepancy function associated with (5.9) is

$$D(P, Q) \equiv |E_P(T) - E_Q(T)|^2 .$$

In particular, with $S$ as in (5.12), we obtain

$$D(P, Q) \equiv |\phi_P(u) - \phi_Q(u)|^2 . \tag{6.7}$$

□

Now note that, if (5.8) holds, we have, in obvious notation,

$$D(P, Q) \equiv \int D_\lambda(P, Q) dW(\lambda). \tag{6.8}$$

**Example 4** We thus have here

$$D(P, Q) \equiv \int |\phi_P(u) - \phi_Q(u)|^2 \, dW(u), \tag{6.9}$$

which will be strict if $W$ has full support. Again, (6.9) is translation invariant.

The design criterion in this example is thus expressible as:

$$\textit{Minimise the expected discrepancy:} \int E\left\{\left|\phi_Y^Z(u) - \phi_\Theta^Z(u)\right|^2\right\} dW(u).$$

□

**Example 5** More generally, in this example we obtain the discrepancy function

$$D(P, Q) \equiv ||P - Q||_K^2, \tag{6.10}$$

which is strict if $K$ is positive definite. This discrepancy function was introduced by Eaton *et al.* (1996).

The design criterion can now be expressed as:

*Minimise the expected squared distance between the sampling distribution and the predictive distribution of $Z$:*

$$E\left\{\left\|P_Y^Z - P_\Theta^Z\right\|_K^2\right\}.$$

□

**Example 6** This is very similar to Example 5, on replacing $\|\cdot\|_K^2$ by $\|\cdot\|_\mu^2$. □

Note that quadratic criteria, as in Examples 5 and 6, lead to discrepancy functions that are symmetric in $P$ and $Q$. An extension of the argument of Savage (1971, §5.2) would suggest that only quadratic-type discrepancies can be both coherent and symmetric.

**Example 7** The discrepancy function $D(P, Q)$ generated by (5.23) is given, for $P, Q \in \mathcal{P}$, by (5.24). □

**Example 8** For this special case of Example 7, we obtain, for $P$ and $Q$ possessing second moments,

$$\begin{aligned} D(P,Q) \equiv & \frac{1}{2}\left\{(\mu_P - \mu_Q)^{\mathrm{T}}\Sigma_Q^{-1}(\mu_P - \mu_Q) + \operatorname{tr}\left(\Sigma_Q^{-1}\Sigma_P\right)\right. \\ & \left. - \log\det\left(\Sigma_Q^{-1}\Sigma_P\right) - p\right\}. \end{aligned} \tag{6.11}$$

This is not strict, vanishing if $\mu_Q = \mu_P$ and $\Sigma_Q = \Sigma_P$.

Again, the design criterion can be re-expressed as the minimisation of $E\{D(P^Z_\Theta, P^Z_Y)\}$. □

# 7 Dependence functions

Let $L(z, a)$ be a loss function, for $z \in \mathcal{Z}$, $a \in \mathcal{A}$, with domain $\mathcal{P}$. Let the space $\mathcal{U}$ be arbitrary, and consider a pair $(Z, U)$ of random variables, with joint distribution $P^{Z,U}$ over $\mathcal{Z}\times\mathcal{U}$. We suppose that, almost surely, the conditional distribution $P^Z_U$ belongs to $\mathcal{P}$. Then the Bayes loss if an action must be chosen immediately is $H(P^Z) = L\left\{P^Z, a(P^Z)\right\}$, while the expected Bayes loss if $U$ can be observed first is $E\left\{H(P^Z_U)\right\}$, the expectation being over the marginal distribution $P^U$ of $U$. Since $H$ is concave and $E\left(P^Z_U\right) = P^Z$, we have $E\left\{H(P^Z_U\right\} \leq H(P^Z)$, with equality if $Z$ and $U$ are independent. The quantity

$$C(P^{Z,U}) := H(P^Z) - E\left\{H\left(P^Z_U\right)\right\} \tag{7.1}$$

may be considered a measure of the dependence of $Z$ on $U$ in $P^{Z,U}$. We call $C$ the *coherent dependence function* corresponding to $L$, and say that $C$ is *strict* if $L$ is. Then $C$ is nonnegative, and vanishes if $Z \perp\!\!\!\perp U$ (and only in this case when $C$ is strict). Note that, even if $\mathcal{U} = \mathcal{Z}$, $C(P^{Z,U})$ need not be symmetric as between $Z$ and $U$. In terms of the original loss function,

$$C(P^{Z,U}) \equiv E_{P^{Z,U}}\left[L\left\{Z, a(P^Z)\right\} - L\left\{Z, a(P^Z_U)\right\}\right] \tag{7.2}$$

defines the "expected value of sample information" contained in $U$ (Raiffa and Schlaifer, 1961, §4.5.2). More generally, $C$ might be defined directly in terms of (7.1), for a given coherent uncertainty function $H$. Note that $C$ is defined for all suitable joint distributions over $\mathcal{Z}\times\mathcal{U}$, for arbitrary sample space $\mathcal{U}$. If $V$ is a one-one function of $U$, $C(P^{Z,V}) = C(P^{Z,U})$. However, $C$ is typically not preserved under transformation of $Z$.

If we start with a proper scoring rule $S$ we obtain

$$C(P^{Z,U}) \equiv E\left\{S(Z, P^Z) - S(Z, P^Z_U)\right\}. \tag{7.3}$$

Strongly equivalent versions of $H$ (or of $L$, or of $S$) lead to the identical $C$. Similarly, starting from a coherent discrepancy function $D$ we have

$$C(P^{Z,U}) \equiv E\left\{D(P_U^Z, P^Z)\right\}. \tag{7.4}$$

Indeed, we have the more general result:

**Lemma 7.1** *For any $Q \in \mathcal{P}$,*

$$C(P^{Z,U}) \equiv E\left\{D(P_U^Z, Q)\right\} - D(P^Z, Q). \tag{7.5}$$

**Proof** Follows from $S(P^Z, Q) \equiv E\left\{S(P_U^Z, Q)\right\}$, and $D(P,Q) \equiv S(P,Q) - H(P)$. □

## 7.1 Predictive design using dependence functions

Returning to our design problem, it follows from (7.1) and the fact that $P^Z$ does not depend on the experiment $\xi$ that the criterion for choice of experiment, minimisation of (4.1), can be alternatively re-expressed as:

*Maximise the dependence $C(P^{Z,Y})$ between*
*the predictand $Z$ and the observable $Y$.*

**Examples 1 and 3** Here we obtain

$$C(P^{Z,Y}) \equiv \operatorname{var}\left\{E(T \mid Y)\right\}.$$

The design criterion is thus:

*Maximise the variance of the predictive expectation of $T$.*

□

**Example 2** Example 2 yields

$$C(P^{Z,Y}) \equiv E\left\{\log p(Z,Y)/p(Z)p(Y)\right\}, \tag{7.6}$$

where $p(z,y)$ denotes the density of $P^{Z,Y}$ at $(Z,Y)$, *etc.* This is the "mutual information" or "cross-entropy" between $Z$ and $Y$. Formula (7.6) treats $Z$ and $Y$ symmetrically and it seems plausible that this is the only case in which this holds.

The design criterion in this case thus becomes:

*Maximise the cross-entropy between $Z$ and $Y$.*

This criterion for experimental design was introduced by Lindley (1956). □

**Example 4** Here we obtain

$$\begin{aligned} C(P^{Z,Y}) &\equiv \int \left[E\left\{|\phi_Y^Z(u)|^2\right\} - |\phi^Z(u)|^2\right] dW(u) \\ &\equiv \int \text{var}\left\{\phi_Y^Z(u)\right\} dW(u). \end{aligned} \tag{7.7}$$

The design criterion is thus to maximise (7.7). □

**Examples 5 and 6** More generally, for Example 5,

$$C(P^{Z,Y}) \equiv E\left\{\|P_Y^Z - P^Z\|_K^2\right\},$$

a generalised variance of the predictive distribution. A similar form holds for Example 6. A good design is thus one leading to predictive distributions which are highly variable, *i.e.* highly sensitive to the data. □

**Examples 7 and 8** In Example 7 we obtain

$$C(P^{Z,Y}) \equiv b(\lambda^Z) - \left(\lambda^Z\right)^{\mathrm{T}} \tau^Z - E\left\{b(\lambda_Y^Z) - \left(\lambda_Y^Z\right)^{\mathrm{T}} \tau_Y^Z\right\},$$

where $\lambda_y^Z$ denotes $\lambda_{P_y^Z}$, *etc.*

The special case of Example 8 gives as design criterion:

$$\textit{Maximise } E\left[\log\left\{\det \Sigma^Z / \det \Sigma_Y^Z\right\}\right].$$

In particular, if the joint distribution of $Y$ and $Z$ is normal in every experiment $\xi$, this becomes:

$$\textit{Maximise } E\left\{\log\det \Sigma^Y\right\}.$$

□

# 8 Prediction and estimation

We have already pointed out that our framework for prediction includes the simpler problem of estimation as a special case. But conversely, any predictive decision problem may itself be always be recast as one of estimation. Thus, given a predictive terminal decision problem, with loss function $L(z, a)$, and a family of distributions $\{P_\theta^Z\} \subseteq \mathcal{P}$ for $Z$ given $\Theta$, we may define an "estimative loss function" $L^\dagger(\theta, a)$ by:

$$L^\dagger(\theta, a) := L(P_\theta^Z, a),$$

*viz.* $E\{L(Z, a)\}$ under $P_\theta^Z$.

If we now truncate Figure 2 at $\nu_4$, and attach there the loss $L^\dagger(\theta, a)$, in accordance with step $(i)$ of §4, this will not further affect steps $(ii)$–$(v)$, and thus

$$
\begin{array}{ccccccccccc}
 & \xi & & y & & Q & & \theta & & & \\
\square & \rightarrow & \bigcirc & \rightarrow & \square & \rightarrow & \bigcirc & \rightarrow & \bullet & & D(P^Z_\theta, Q) \\
\nu_0 & & \nu_1 & & \nu_2 & & \nu_3 & & \nu_4 & &
\end{array}
$$

Figure 4: Estimative design

we shall obtain the same optimal $\xi$ and $a$. Only now the problem has been cast as one of estimation, with loss function $L^\dagger(\theta, a)$.

If we analyse this new estimation problem as a decision problem in its own right, we find, in an obvious notation, the following associated functions (where, for $\Pi$ a distribution for $\Theta$, $P^Z_\Pi := \int P^Z_\theta d\Pi(\theta)$ denotes the induced predictive distribution for $Z$ in the original problem):

$$
\begin{aligned}
L^\dagger(\Pi, a) &\equiv L(P^Z_\Pi, a) \\
a^\dagger(\Pi) &\equiv a(P^Z_\Pi) \\
S^\dagger(\theta, \Gamma) &\equiv S(P^Z_\theta, P^Z_\Gamma) \\
S^\dagger(\Pi, \Gamma) &\equiv S(P^Z_\Pi, P^Z_\Gamma) \\
H^\dagger(\Pi) &\equiv H(P^Z_\Pi) \\
D^\dagger(\Pi, \Gamma) &\equiv D(P^Z_\Pi, P^Z_\Gamma).
\end{aligned}
$$

Also, so long as $V \perp\!\!\!\perp Z \mid \Theta$,

$$C^\dagger(P^{\Theta,V}) \equiv C(P^{Z,V}).$$

Hence the above relations define proper scoring rules and coherent uncertainty functions, discrepancy functions and dependence functions on the parameter space, which could be applied directly in the problem of estimation. An interesting problem is to characterise when, given (say) an estimative uncertainty function $H^\dagger$, it can arise as above from some loss function in a predictive problem with given $\{P^Z_\theta\}$.

Now consider again the analysis of §6.1. The terms involving $H$ in $(i)'' - (iv)''$ do not affect the minimising values for $Q$ in $(iii)''$ and for $\xi$ in $(v)''$. Consequently, the identical solutions are obtained from the decision problem of Figure 4, where we have ignored $H$ altogether, and collapsed out the final stage using $(i)''$.

In this reformulation, the task can be regarded as that of predicting the sampling distribution $P^Z_\theta$, with loss function $L'(\theta, Q)$ given by the discrepancy $D(P^Z_\theta, Q)$ between the true and the quoted distribution. Note that $L'(\theta, Q)$ is an "estimative" loss function, in so much as it depends on the value $\theta$ of the parameter $\Theta$; and that the action space is $\mathcal{P}^Z$, the space of distributions for $Z$. The resulting $S'$, $H'$, $D'$ and $C'$ are strongly equivalent to $S^\dagger$, $H^\dagger$, $D^\dagger$ and $C^\dagger$ above, and of course the optimal $\xi$ and $a$ are the same.

We can thus reformulate the predictive design criterion once again, as:

*Solve the estimative design problem, where the loss function is*

$$L'(\theta, Q) \equiv D(P_\theta^Z, Q).$$

**Example 1** In Example 1, the equivalent estimative loss can be taken as:

$$L^\dagger(\theta, a) \equiv \sigma_\theta^2 + (\mu_\theta - a)^2,$$

or, equivalently, $(\mu_\theta - a)^2$, where $\mu_\theta$, $\sigma_\theta^2$ denote, respectively, the mean and variance of $T \equiv T(Z)$ under $P_\theta^Z$.

Alternatively, with $\mathcal{P}^Z$ as action space, we can use:

$$L'(\theta, Q) \equiv \{\mu_\theta - E_Q(T)\}^2.$$

□

**Examples 2–6** In Example 2 we can use, as estimative loss, the Kullback-Liebler distance of $Q$ from $P_\theta^Z$. For Example 5, this is replaced by the quadratic distance, $\|P_\theta^Z - Q\|_K^2$, or $\|P_\theta^Z - Q\|_\mu^2$ in Example 6. Examples 3 and 4 are special cases of this, the latter giving estimative loss

$$L'(\theta, Q) \equiv \int \left|\phi_\theta^Z(u) - \phi_Q(u)\right|^2 dW(u).$$

This loss function is the basis of the approach taken by Eaton *et al.* (1996). □

**Examples 7 and 8** For Example 7 the estimative loss will depend on $\theta$ only through the expectation $\tau_\theta$ of $T \equiv t(Z)$ under $P_\theta^Z$. In the particular case of Example 8 we obtain:

$$\begin{aligned} L'(\theta, Q) \equiv \; & \frac{1}{2}\Big\{(\mu_\theta - \mu_Q)^{\mathrm{T}}\Sigma_Q^{-1}(\mu_\theta - \mu_Q) + \operatorname{tr}\left(\Sigma_Q^{-1}\Sigma_\theta\right) \\ & - \log\det\left(\Sigma_Q^{-1}\Sigma_\theta\right) - p\Big\}, \end{aligned} \tag{8.1}$$

where $\mu_\theta$ and $\Sigma_\theta$ denote, respectively, the mean vector and dispersion matrix of the distribution of $Z$ given $\theta$. □

# 9 Characterisations of coherence

So far we have defined coherent uncertainty functions, discrepancy functions and dependence functions simply as those which can arise, according to the definitions introduced, starting from some terminal decision problem with its loss function $L$. This raises the general question: Given such a function (on a convex domain $\mathcal{P}$) which is purported to be coherent, how may one check this? This question is important if, say, some one has proposed a design criterion of the same general form of those considered above, based on some sort of uncertainty function, discrepancy function, *etc.* If we can show that the function used is coherent, we shall then know that the proposed design criterion is coherent, in the sense of being equivalent to conducting a full Bayesian decision analysis based on some loss function. Such an investigation of coherence has been carried out in detail by Dawid and Sebastiani (1996) for some commonly recommended criteria for optimal design of regression experiments.

We would thus wish to characterise necessary and/or sufficient conditions which can be used to confirm or refute coherence. We therefore now turn to a closer consideration of the characteristics of coherence. While necessary conditions for coherence may be expressed fairly simply, sufficient conditions are generally more delicate, and we shall only characterise these in somewhat limited mathematical contexts. However, there is never any harm in attempting to apply these sufficient conditions heuristically in other circumstances, since they can often suggest a suitable loss function or proper scoring rule to examine as a possible source of the function in question, and this may then be checked directly. Of course, once one has exhibited such a loss function, coherence is established.

**Definitions and assumptions** Let $\mathcal{P}$ be a convex family of distributions over $\mathcal{Z}$. We shall confine attention to distributions $P \in \mathcal{P}$. If $\mathcal{V}$ is an affine space, we shall call a function $F : \mathcal{P} \rightarrow \mathcal{V}$ *affine* if, for $\tilde{P}$ a random element of $\mathcal{P}$, $E\{F(\tilde{P})\} = F\{E(\tilde{P})\}$.

Suppose that, for a loss function $L$ with domain $\mathcal{P}$, we define, for any $a \in \mathcal{A}$, a function $F_a : \mathcal{P} \rightarrow \mathbb{R}$ by:

$$F_a(P) \equiv L(P, a), \tag{9.1}$$

as given by (3.1). Then for each $a$ we have an affine function $F_a : \mathcal{P} \rightarrow \mathbb{R}$, and the coherent Bayesian solution to the decision problem when $Z \sim P \in \mathcal{P}$ is to choose $a$ to minimise $F_a(P)$. This last sentence might be taken, in a suitable axiomatic approach which limits attention to $\mathcal{P}$, as the very definition of coherence; it is certainly enough to stand as a basis for all the analysis of this paper. It differs from the usual approach simply because the functions $F_a$ might not be *presentable*, defined as follows:

**Definition 9.1** *We call $F : \mathcal{P} \rightarrow \mathcal{V}$* presentable *if there exists a function $f : \mathcal{Z} \rightarrow \mathcal{V}$ such that $F(P) \equiv E_P\{f(Z)\}$.*

Every presentable function is affine, but the validity of the converse in general is a matter for functional-analytical examination of the specific spaces and functions involved. Presentability of the "expected loss" functions $F_a$ implies the existence of a loss function $L$ such that (9.1) holds, and so allows us to phrase the coherence requirement in the usual way, in terms of actual expectations of the loss. If this is regarded as important, we generally need to impose additional conditions.

The following simple lemma establishes this converse in the special case that the family $\mathcal{P}$ contains all one-point distributions.

**Lemma 9.1** *Let $\mathcal{P}$ be the family of all probability distributions over $\mathcal{Z}$, and let $F : \mathcal{P} \to \mathcal{V}$. Then if $F$ is affine, $F$ is presentable.*

**Proof** For $z \in \mathcal{Z}$, $\mathcal{P}$ contains $\delta_z$, the point mass at $z$. Define $f : \mathcal{Z} \to \mathcal{V}$ by: $f(z) \equiv F(\delta_z)$. For $P \in \mathcal{P}$, consider the random element $\tilde{P}$ of $\mathcal{P}$ given by: $\tilde{P} = \delta_Z$, where $Z \sim P$. Then $E(\tilde{P}) = P$. So $F(P) = F\{E(\tilde{P})\} = E\{F(\tilde{P})\}$, since $F$ is affine. Hence $F(P) = E_P\{f(Z)\}$, so that $F$ is presentable. □

We shall not go any more deeply into the validity of such results in more general settings, but henceforth shall simply assume that, for the types of spaces and functions under consideration, we can show that every affine function is presentable. Even when this fails, we can still proceed if we are willing to accept the weaker interpretation of coherence introduced above.

# 10 Properties of uncertainty functions

We have already pointed out that a coherent uncertainty function $H : \mathcal{P} \to \mathbb{R}$ must be concave In fact more is true. Coherence of $H$ is equivalent to its expressibility as $H(P) \equiv S(P,P)$ for some proper scoring rule $S$. For fixed $Q \in \mathcal{P}$, define $\phi_Q : \mathcal{P} \to \mathbb{R}$ by: $\phi_Q(P) \equiv S(P,Q)$. Then we have:

$$\phi_Q \text{ is affine, and } H(P) \leq \phi_Q(P), \text{ all } P \in \mathcal{P}, \text{ with equality at } P = Q. \quad (10.1)$$

Consider first the case that $\mathcal{Z}$ is finite with $N$ values, with $\mathcal{P}$ the set of all distributions over $\mathcal{Z}$. Then we can represent $\mathcal{P}$ as $\mathcal{S}$, the unit simplex in $\mathbb{R}^N$, and the affine function $\phi_Q$ of (10.1) by:

$$\phi_Q(P) \equiv H(Q) + c_Q^{\mathrm{T}}(P - Q)$$

for some $c_Q \in \mathbb{R}^N$, which we may take to lie in $\mathcal{M} := \{x : \sum x_i = 0\}$. Thus

$$H(P) \leq H(Q) + c_Q^{\mathrm{T}}(P - Q). \quad (10.2)$$

For fixed $Q$, varying $P$, the right-hand side of (10.2) defines a "supporting hyperplane" to $H$ at $Q$, and concavity of $H$ follows easily. The linear function $x \mapsto c_Q^{\mathrm{T}} x$ on $\mathcal{M}$ is a "super-gradient" to $H$ at $Q$.

Conversely, any concave $H$ has the supporting hyperplane property (10.2) at each $Q$ in the interior of $\mathcal{S}$, for suitable super-gradient $c_Q \in \mathcal{M}$, not necessarily unique. When $H$ is differentiable over $\mathcal{S}$ at $Q$, we have $c_Q = \nabla H(Q)$ (uniquely defined in $\mathcal{M}$).

In more general spaces, (10.1) again implies concavity, as expressed in the property: if $\tilde{P}$ is a random element over $\mathcal{P}$, then

$$E\{H(\tilde{P})\} \leq H\{E(\tilde{P})\}. \tag{10.3}$$

The converse need not hold: see Hendrickson and Buehler (1971), Example 4.1, for a counter-example. However, in suitable spaces, and under suitable continuity conditions, it will do so (see, *e.g.*, *ibid.*, Theorem 4.1; Johnson, 1991). Under differentiability, we shall have $\phi_Q(P) \equiv H(Q) + \psi_Q(P - Q)$, with $\psi_Q$ the Gâteaux derivative of $H$ at $Q$, and the "supporting hyperplane" is just the "tangent hyperplane".

Suppose now that $H$ is concave on $\mathcal{P}$, and that we can further show (either directly, or by appealling to suitable general conditions for this for the case at hand) the supporting hyperplane property (10.1). Since we are assuming that $\phi_Q$, being affine, is therefore presentable, there then exists $f_Q : \mathcal{Z} \rightarrow \mathbb{R}$ such that

$$\phi_Q(P) = E_P\{f_Q(Z)\}. \tag{10.4}$$

Defining

$$S(z, Q) \equiv f_Q(z), \tag{10.5}$$

and, as usual, $S(P, Q) \equiv E_P\{S(Z, Q)\} \equiv \phi_Q(P)$, we thus have

$$H(P) \leq S(P, Q), \text{ all } P, Q; \tag{10.6}$$

$$\text{and } H(P) \equiv S(P, P). \tag{10.7}$$

It follows that $S$ is a proper scoring rule giving rise to $H$, so that $H$ is coherent.

When the above analysis applies we can further express the associated discrepancy function $D$ in terms of $H$:

$$D(P, Q) \equiv \phi_Q(P) - H(P). \tag{10.8}$$

In summary, we see that, under appropriate further assumptions and conditions, any concave real function on $\mathcal{P}$ can be used as a coherent uncertainty function.

### Examples

**Example 1** Suppose that we wish to verify that the uncertainty function

$$H(P) := \text{var}_P(T)$$

(with domain $\mathcal{P} = \{P : \text{var}_P(T) < \infty\}$) is coherent. It is readily seen that this is concave in $P \in \mathcal{P}$.

We have

$$\mathrm{var}_P(T) = E_P\left[\{T - E_Q(T)\}^2\right] - \{E_P(T) - E_Q(T)\}^2.$$

Thus

$$\mathrm{var}_P(T) \leq \phi_Q(P) := E_P\left[\{T - E_Q(T)\}^2\right], \text{ all } P, Q.$$

Since the right-hand side is affine in $P$, and we have equality when $P = Q$, we have shown that (10.1) holds. Moreover, $\phi_Q$ is clearly presentable, and we recover $S(z, Q) \equiv \{t(z) - E_Q(T)\}^2$ as a proper scoring rule giving rise to the uncertainty function $\mathrm{var}_P(T)$, which is thus shown to be coherent. □

Example 5 (and thus its special cases, Examples 3 and 4, and variant, Example 6) may be treated in similar fashion.

**Example 2** Now we investigate the coherence of the uncertainty function $H(P)$ given by the entropy, (3.13), on the (convex) domain where this is finite. Let $p(z) \equiv q(z) + \epsilon r(z)$, where $\int r(z) d\mu = 0$. Proceeding informally, we obtain

$$\begin{aligned} H(P) &= H(Q) - \epsilon \int \{1 + \log q(z)\} r(z) d\mu + o(\epsilon) \\ &= H(Q) - \epsilon \int \{\log q(z)\} r(z) d\mu + o(\epsilon). \end{aligned}$$

Hence the "directional derivative" at $Q$ in the direction $r$ is $-\int\{\log q(z)\} r(z) d\mu$, and we are therefore led to examine the "tangent hyperplane" to $H$ at $Q$:

$$\begin{aligned} \phi_Q(P) &\equiv H(Q) - \int \{\log q(z)\}\{p(z) - q(z)\} d\mu \\ &= -E_P\{\log q(Z)\} \end{aligned}$$

and, correspondingly,

$$S(z, Q) \equiv -\log q(z).$$

It is then an easy matter to verify (as indeed has already been done) that $S$ is a proper scoring rule, and generates the entropy as its corresponding uncertainty function, which must therefore be coherent. Any liberties taken in the above "derivation" are then entirely validated by this verification.

This argument shows, in particular, the coherence of the design criterion based on minimising the expected entropy of the predictive (or posterior) distribution. □

**Example 7** Let $b$ be a closed, convex and differentiable function on the open convex set $\mathcal{L} \subseteq \mathbb{R}^k$, with differential $\beta := \nabla b$. (Similar results are obtainable if $b$ is merely subdifferentiable. See *e.g.* Barndorff-Nielsen, 1978, Chapter 5, for relevant definitions and background on convex functions). Let $\mathcal{T} \subseteq \mathbb{R}^k$ be the image of $\mathcal{L}$ under $\beta$. Then $\mathcal{T}$ is open and convex, and $\beta : \mathcal{L} \to \mathcal{T}$ is invertible. The *conjugate function* $b^* : \mathcal{T} \to \mathcal{L}$ of $b$ is defined by:

$$b^*(\tau) \equiv \sup_{\delta \in \mathcal{L}} \left\{ \delta^{\mathrm{T}} \tau - b(\delta) \right\}, \tag{10.9}$$

and the supremum in (10.9) is attained at $\delta = \lambda$, where $\beta(\lambda) = \tau$. Then $b^*$ is a closed convex function. From (10.9) we deduce that

$$b^*(\gamma) \geq b^*(\tau) + \lambda^{\mathrm{T}}(\gamma - \tau) \tag{10.10}$$

for all $\gamma \in \mathcal{T}$.

Now let $T \equiv t(Z)$, where $t : \mathcal{Z} \to \overline{\mathcal{T}}$, $\overline{\mathcal{T}}$ being the closure of $\mathcal{T}$, and take as domain $\mathcal{P}$ the set of distributions $P$ for $Z$ such that $\tau_P := E_P(T) \in \mathcal{T}$. Define $H$ on $\mathcal{P}$ by

$$\begin{aligned} H(P) &\equiv -b^*(\tau_P) \\ &\equiv b(\lambda_P) - \lambda_P^{\mathrm{T}} \tau_P, \end{aligned} \tag{10.11}$$

with $\beta(\lambda_P) = \tau_P$ (*cf*. (5.25)). Then $H$ is concave on $\mathcal{P}$. From (10.10) and (10.11), we find

$$H(P) \leq H(Q) - \lambda_Q^{\mathrm{T}}(\tau_P - \tau_Q). \tag{10.12}$$

The right-hand side of (10.12) is affine in $P$, so that (10.1) is satisfied. Clearly it is also presentable, and thus we have demonstrated that $H$ given by (10.11) is a coherent uncertainty function. The associated proper scoring rule is

$$\begin{aligned} S(z, Q) &\equiv H(Q) - \lambda_Q^{\mathrm{T}} \{t(z) - \tau_Q\} \\ &\equiv b(\lambda_Q) - \lambda_Q^{\mathrm{T}} t(z) \end{aligned}$$

(*cf*. (5.23)). □

**Example 8** We now investigate the coherence of the uncertainty function

$$H(P) \equiv \log \det \Sigma_P, \tag{10.13}$$

equivalent to (5.28), defined for the set $\mathcal{P}$ of distributions $P$ over $\mathcal{Z} = \mathbb{R}^p$ having finite mean $\mu_P$ and dispersion matrix $\Sigma_P$. Although this can be regarded as a special case of Example 7, we shall analyse it from first principles.

We have, for $\alpha, \beta \geq 0$, $\alpha + \beta = 1$,

$$\begin{aligned} \Sigma_{\alpha P + \beta Q} & \quad \alpha \Sigma_P + \beta \Sigma_Q + \alpha\beta(\mu_P - \mu_Q)(\mu_P - \mu_Q)^{\mathrm{T}} \\ &\geq \alpha \Sigma_P + \beta \Sigma_Q \end{aligned}$$

in the Loewner ordering. Since the function $\Sigma \to \log \det \Sigma$ is concave (see *e.g.* Dawid and Sebastiani, 1996), $H$ is concave in $P$. And, since $H$ can be regarded as a concave function of $(E_P(Z), E_P(ZZ^{\mathrm{T}}))$, which ranges over a finite-dimensional open convex set, $H$ must satisfy the supporting hyperplane property (10.1), and hence be coherent. It only remains to identify the function $\phi_Q$, and hence the associated proper scoring rule $S$.

Let $R$ be a signed measure of zero total mass such that $P := Q + R \in \mathcal{P}$, and consider $P_\epsilon := Q + \epsilon R = (1-\epsilon)Q + \epsilon P$. Then

$$\begin{aligned} \Sigma_{P_\epsilon} &= (1-\epsilon)\Sigma_Q + \epsilon\Sigma_P + \epsilon(1-\epsilon)(\mu_P - \mu_Q)(\mu_P - \mu_Q)^{\mathrm{T}} \\ &= \Sigma_Q + \epsilon\left\{(\Sigma_P - \Sigma_Q) + (\mu_P - \mu_Q)(\mu_P - \mu_Q)^{\mathrm{T}}\right\} + o(\epsilon). \qquad (10.14) \end{aligned}$$

Now the directional derivative of the matrix function $\log\det\Sigma$ in the direction $\Gamma$, where $\Sigma$ and $\Sigma + \Gamma$ are both positive definite symmetric matrices, is $\mathrm{tr}\Sigma^{-1}\Gamma$, as may be deduced from Sebastiani (1995). That is, for small $\epsilon$,

$$\log\det(\Sigma + \epsilon\Gamma) = \log\det(\Sigma) + \epsilon\mathrm{tr}\Sigma^{-1}\Gamma + o(\epsilon).$$

Inserting (10.14) thus gives

$$\log\det\Sigma_{P_\epsilon} = \log\det\Sigma_Q + \epsilon\mathrm{tr}\left[\Sigma_Q^{-1}\left\{(\Sigma_P - \Sigma_Q) + (\mu_P - \mu_Q)(\mu_P - \mu_Q)^{\mathrm{T}}\right\}\right] + o(\epsilon). \qquad (10.15)$$

The affine approximation to $\log\det\Sigma_P$ around $Q$ in the direction of $R$ is then given by ignoring the $o(\epsilon)$ term in (10.15). In particular, taking $\epsilon = 1$ yields

$$\log\det\Sigma_Q + \mathrm{tr}\left\{\Sigma_Q^{-1}E_P(Z - \mu_Q)(Z - \mu_Q)^{\mathrm{T}}\right\} - p$$

as the approximation at $P$. Since this is fully affine (and presentable) in $P$, it defines the tangent hyperplane, and hence the supporting hyperplane $\phi_Q(P)$ in (10.1). The corresponding scoring rule is

$$S(z, Q) \equiv \log\det\Sigma_Q + \mathrm{tr}\left\{\Sigma_Q^{-1}(z - \mu_Q)(z - \mu_Q)^{\mathrm{T}}\right\} - p,$$

equivalent to (5.27), which has already been shown to be proper. Hence the coherence of the "$D$-optimality" uncertainty function $\log\det\Sigma_P$ has been established. □

The uncertainty function (10.13) is just one of a wide variety of coherent uncertainty functions which depend only on the dispersion matrix of $Z$; these are characterised by Dawid and Sebastiani (1996).

# 11 Properties of discrepancy functions

Let $D$ be a coherent discrepancy function over $\mathcal{P}$. Then for $Q, Q_0,\ P \in \mathcal{P}$,

$$D(P, Q) - D(P, Q_0) \equiv \int \{S(z, Q) - S(z, Q_0)\}\, dP(z).$$

In particular we see that, for any fixed $Q, Q_0 \in \mathcal{P}$,

$$D(P,Q) - D(P,Q_0) \text{ is affine in } P. \tag{11.1}$$

Conversely, suppose $D$ is a discrepancy function over $\mathcal{P}$, *i.e.* for $P, Q \in \mathcal{P}$, $D(P,Q) \geq 0$ with $D(P,P) \equiv 0$; and suppose further that (11.1) holds. Since we are assuming that any affine function is presentable, there exists a function $f(z; Q, Q_0)$ such that

$$D(P,Q) - D(P,Q_0) \equiv E_P \{f(Z; Q, Q_0)\}. \tag{11.2}$$

Fix $Q_0 \in \mathcal{P}$, and define

$$S(z; Q) \equiv f(z; Q, Q_0). \tag{11.3}$$

So

$$\begin{aligned} S(P,Q) &\equiv D(P,Q) - D(P,Q_0) \\ S(P,P) &\equiv -D(P,Q_0). \end{aligned}$$

Then $S(P,Q) - S(P,P) \equiv D(P,Q) \geq 0$. Thus $S$ is a proper scoring rule, and $D$ is its corresponding discrepancy function, which is thereby shown to be coherent. The associated uncertainty function is

$$H(P) \equiv -D(P,Q_0). \tag{11.4}$$

Different choices for $Q_0$ yield strongly equivalent problems.

**Example 2** Take $D \equiv K$, the Kullback-Liebler distance (3.10), where $\mu$ is an arbitrary measure. Fix a probability measure $Q_0 \ll \mu$, and take as domain $\mathcal{P} = \{\mathcal{P} : \mathcal{P} \ll \mathcal{Q}_{\prime}, \mathcal{K}(\mathcal{P}, \mathcal{Q}_{\prime})\infty\}$. We find

$$K(P,Q) - K(P,Q_0) \equiv E_P [\log \{q_0(Z)/q(Z)\}],$$

which is affine, and clearly presentable, in $P$. Hence $K$ is a coherent discrepancy function. The above analysis then produces the associated proper scoring rule $S(z,Q) \equiv \log q_0(z) - \log q(z)$, strongly equivalent to $-\log q(z)$; and coherent uncertainty function $H(P) \equiv -K(P,Q_0)$, the entropy with respect to $Q_0$. □

Again, the other examples may be confirmed similarly.

## 11.1 A-coherence

Aitchison (1990) introduced another coherence property for discrepancy functions, essentially equivalent to the following.

**Definition 11.1** The discrepancy measure $D$ is called *A-coherent* if, for $\tilde{P}$ a random distribution over $\mathcal{P}$, $E\left\{D(\tilde{P}, Q)\right\}$ is minimised in $Q$ when $Q = E(\tilde{P})$.

This is in fact closely related to our definition of coherence. Thus we have:

**Theorem 11.1** *If $D$ is coherent, then $D$ is A-coherent.*

**Proof** Fix $Q_0 \in \mathcal{P}$, and define $U(P, Q) \equiv D(P, Q) - D(P, Q_0)$. By coherence of $D$, $E\left\{U(\tilde{P}, Q)\right\} \equiv U\left\{E(\tilde{P}), Q\right\}$. Thus $E\left\{D(\tilde{P}, Q)\right\} \equiv D\left\{E(\tilde{P}), Q\right\} + \left[E\left\{D(\tilde{P}, Q_0)\right\} - D\left\{E(\tilde{P}), Q_0\right\}\right]$, and this is minimised for $Q = E(\tilde{P})$. □

Since there is a wide variety of coherent discrepancy functions, Theorem 11.1 refutes Aitchison's conjecture that only the Kullback-Liebler distance has the A-coherence property. Eaton *et al.* (1995) have noted this property for the coherent discrepancy (6.9).

The converse to Theorem 11.1 can be proved under suitable conditions of regularity and differentiability: see Dawid (1994) for a simple version.

# 12 Properties of dependence functions

We here note some properties of coherent dependence functions, of general interest.

**Lemma 12.1** *Let $(Z, U, V)$ have a joint distribution such that $Z \perp\!\!\!\perp V \mid U$. Then $C(P^{Z,U}) \geq C(P^{Z,V})$.*

**Proof** Since $E(P^Z_U \mid V) = P^Z_V$, and $H$ is concave, $E\left\{H(P^Z_U) \mid V\right\} \leq H(P^Z_V)$, so that $E\left\{H(P^Z_U)\right\} \leq E\left\{H(P^Z_V)\right\}$, and the result follows from (7.1). □

**Corollary 12.1** *If $V$ is a function of $U$, $C(P^{Z,U}) \geq C(P^{Z,V})$.*

**Corollary 12.2** *If $H$ in (7.1) is strict, $Z \perp\!\!\!\perp V \mid U$ and $C(P^{Z,U}) = C(P^{Z,V})$, then $Z \perp\!\!\!\perp U \mid V$.*

**Proof** In this case $E\left\{H(P^Z_{U,V}) \mid V\right\} = H(P^Z_V)$ almost surely, so by strictness $P^Z_{U,V} = P^Z_V$. □

Suppose now that $C$ is bounded, with domain the set of all distributions for $Z$. Then we have

**Corollary 12.3** *For fixed $P^Z$, $C(P^{Z,U})$ is maximised for $U = Z$.*

| | $S$ | $H$ | $D$ | $C$ |
|---|---|---|---|---|
| $L$ | (5.3) | (3.3) | (6.3) | (7.2) |
| $S$ | | (5.2) | (6.1) | (7.3) |
| $H$ | (10.1)<br>(10.4)<br>(10.5) | | (10.1)<br>(10.8) | (7.1) |
| $D$ | (11.2)<br>(11.3) | (11.4) | | (7.4) |
| $C$ | | (12.1) | | |

Table 1: *Equations specifying the relationships between a loss function $L(x, a)$, a proper scoring rule $S(x, P)$, a proper uncertainty function $H(P)$, a proper discrepancy function $D(P, Q)$, and a proper dependence function $C(P^{X,U})$.*

None of these properties leads easily to a characterisation of coherence for a dependence function. However, in the setting of Corollary 12.3, we obtain from (7.1)

$$H^*(P^Z) \equiv C(P^{Z,Z}), \tag{12.1}$$

where $H^*$ is the canonical version of $H$, given by (3.8). Hence, given $C$, $H$ may be recovered up to equivalence. It is then in principle possible to check that (12.1) yields a coherent uncertainty function, and that (7.1) holds, which together constitute necessary and sufficient conditions for $C$ to be coherent.

For cases such as the cross-entropy, where $C(P^{Z,Z})$ is undefined, such a routine is not available, and we do not currently have a simple characterisation of coherence for general dependence functions.

# 13 Summary of relationships

Table 1 summarises the major relationships between the various functions discussed in this paper. Each entry gives the relevant equation(s) relating the quantity labelling its column to that labelling its row, such a relationship being determined up to strict equivalence. Given any of the other quantities, a suitable definition of a related loss function $L$ may be taken to be the same as $S$. Expressions for $S$ and $D$ in terms of $C$ are not given explicitly, but may be obtained by first expressing them in terms of $H$.

# 14 Conclusions

The coherent approach to decision analysis is so fundamental and versatile that its implications deserve to be taken extremely seriously. We have seen how, using this approach, it is possible to define specific kinds of measures of uncertainty, discrepancy and dependence, which we have likewise termed coherent, and which have special properties which we have characterised. In particular, in Example 2 we have seen how one particular decision problem, based on the logarithmic score, gives rise to a number of fundamental concepts: entropy; Kullback-Liebler distance; and cross-entropy. Since the logarithmic score is essentially the same as log-likelihood, which is one of the most fundamental concepts of Theoretical Statistics, the fundamental nature of these concepts may be regarded as inherited from likelihood.

For the particular problem of Bayesian experimental design, we have seen how we can base our design criterion on any of the derived functions, so long as it is coherent, related functions giving rise to the same solution.

The coherent functions studied here also have applications in quite different areas. For example, Dawid (1991, Section 8) uses coherent discrepancy functions in a fundamental way to analyse the performance of a sequential forecasting system based on an incorrect model. Also, other types of coherent function may be introduced. Thus, by considering neighbouring distributions, we can associate, with any discrepancy function, a metric over $\mathcal{P}$, considered as a differential-geometric manifold (see, *e.g.,* Amari, 1985). *Coherent* metrics, arising from coherent discrepancies, again have special properties, which may be characterised, and each such coherent metric determines its associated decision problem up to equivalence. In particular, the Kullback-Liebler distance (Example 2) generates the fundamental Fisher information metric. Eguchi (1983) has characterised a wide class of discrepancy functions having this property, but the Kullback-Liebler distance is distinguished among these as the only one which is coherent.

We see then that there is a variety of types of coherent function, each with an enormous potential range of statistical applications. It would thus appear to be worthwhile, in any study involving general functions describing uncertainty, discrepancy, *etc.*, to investigate in depth the consequences of assuming such functions to possess the additional properties implied by coherence.

# Acknowledgments

I have benefited greatly from discussions with P. Sebastiani and A. Giovagnoli.